\documentclass[11pt,a4paper]{article}

\usepackage{natbib}
\usepackage{amsmath,amsthm,amssymb, thm-restate}
\usepackage{todonotes}
\usepackage{physics}
\usepackage{float}
\usepackage[margin=1in]{geometry}
\usepackage[boxed,vlined,linesnumbered,ruled,resetcount]{algorithm2e}
\usepackage{stmaryrd}
\usepackage{systeme}
\usepackage{authblk}
\usepackage{tocbibind}
\usepackage[toc,page]{appendix}
\usepackage{bbold}
\usepackage{bm}
\usepackage{comment}
\usepackage{rotating}
\usepackage{pdflscape}
\usepackage{physics}
\usepackage{mathtools}
\usepackage{soul}

\usepackage{graphicx}%
\usepackage{multirow}%
\usepackage{amsmath,amssymb,amsfonts, stmaryrd}%
\usepackage{amsthm}%
\usepackage{thmtools}
\usepackage{thm-restate}
\usepackage{mathrsfs}%
\usepackage{xcolor}%
\usepackage{textcomp}%
\usepackage{manyfoot}%
\usepackage{booktabs}%
\usepackage{multirow, multicol, caption, makecell, nicematrix, booktabs}
\usepackage[graphicx]{realboxes}
\usepackage{subfigure}
\usepackage{placeins}

\usepackage{hyperref}
\hypersetup{
colorlinks=true,
breaklinks=true,
urlcolor= blue,
linkcolor= black, 
citecolor=black}

\newcommand{\N}{\mathbb{N}}
\newcommand{\R}{\mathbb{R}}
\newcommand{\Z}{\mathbb{Z}}
\newcommand{\CMax}{\textup{CMax}}
\newcommand{\piEq}{\pi^\text{eq}}
\newcommand{\piInf}{\pi^\text{inf}}
\newcommand{\piSup}{\pi^\text{sup}}
\newcommand{\EnsBinomial}[2]{I_{#2}^{#1}}
\newcommand{\absolute}[1]{\left \lvert #1 \right \rvert}

\newcommand{\ml}[1]{{ \color{blue} #1}}

\newtheorem{theorem}{Theorem}
\numberwithin{theorem}{section}
\newtheorem{definition}[theorem]{Definition}

\newtheorem{proposition}[theorem]{Proposition}

\title{Quadratic versus Polynomial Unconstrained Binary Models for Quantum Optimization illustrated on Railway Timetabling}
\date{}
\author[1,2,3]{\underline{Camille Grange}}
\author[3,4]{Marion Lavignac}
\author[3]{Valentina Pozzoli}
\author[2]{Eric Bourreau}

\affil[1]{\small{LITIS, Normandie Univ, UNIHAVRE, 25 rue Philippe Lebon, Le Havre, France, email: camille.grange@univ-lehavre.fr}}
\affil[2]{\small{LIRMM, Univ Montpellier, CNRS, 161 rue Ada, Montpellier, France, email: eric.bourreau@lirmm.fr}}
\affil[3]{\small{SNCF, Technology, Innovation and Group Projects Department, 1 avenue François Mitterand, Saint-Denis, France, email: valentina.pozzoli@sncf.fr}}
\affil[4]{\small{CentraleSupélec, Paris-Saclay University, 3 rue Joliot-Curie, Gif-sur-Yvette, France, email: marion.lavignac@student-cs.fr}}
\begin{document}

\emergencystretch 3em 

\setlength{\parskip}{0mm}
\setlength{\abovedisplayskip}{3mm}
\setlength{\belowdisplayskip}{3mm}
\setlength{\parindent}{4mm}
\allowdisplaybreaks

\maketitle
\setcounter{page}{1}
\renewcommand{\thepage}{\arabic{page}}

\begin{abstract}
    Quantum Approximate Optimization Algorithm (QAOA) is one of the most short-term promising quantum-classical algorithm to solve unconstrained combinatorial optimization problems. It alternates between the execution of a parametrized quantum circuit and a classical optimization. There are numerous levers for enhancing QAOA performances, such as the choice of quantum circuit meta-parameters or the choice of the classical optimizer. In this paper, we stress on the importance of the input problem formulation by illustrating it with the resolution of an industrial railway timetabling problem. Specifically, 
    we present a generic method to reformulate any polynomial problem into a Polynomial Unconstrained Binary Optimization (PUBO) problem, with a specific formulation imposing penalty terms to take binary values when the constraints are linear. We also provide a generic reformulation into a Quadratic Unconstrained Binary Optimization (QUBO) problem. We then conduct a numerical comparison between the PUBO with binary penalty terms and the QUBO formulations proposed on a railway timetabling problem solved with QAOA. Our results illustrate that the PUBO reformulation outperforms the QUBO one for the problem at hand.
    
    \noindent\textbf{keywords:} Combinatorial optimization, Quantum optimization, QAOA, Railway timetable, Unconstrained binary model
\end{abstract}

\section{Introduction}
\label{sec:intro}

Today, a major part of optimization problems in the industrial context are challenging, with resolutions that remain difficult using current classical methods. The difficulty in solving these problems, namely finding the optimal solution(s), stems both from the inherent combinatorial complexity, often NP-hard, and from the industrial scale of the instances involved. Facing these obstacles, quantum computing is expected to both improve solution quality and reduce computation time. Such expectations are not demonstrated on problems of industrial, practical scale yet, and represent ongoing research. Some operational search methods using quantum algorithms present theoretical advantages compared to classical ones~\citep{Ambainis, GrangeEJOR, Montanaro20_BranchAndBound, Simplex, Kerenidis_LP}, but they require high-quality quantum resources to be implemented. Specifically, they need quantum computers with many qubits that can interact two by two, and quantum operations that can be applied in a row without generating too much noise. However, the current quantum computers do not respect these criteria, encouraging the researchers to look into \emph{lighter} algorithms, waiting for more powerful machines to be built. They consist of hybrid algorithms that take advantage of both quantum and classical resources. Indeed, quantum resources can focus on the critical aspects of the algorithm that utilize quantum information theory, while the more computationally demanding tasks are handled by classical resources.
Such algorithms are metaheuristics, the most famous being the Quantum Approximate Optimization Algorithm (QAOA) introduced by~\cite{QAOAFarhi} to solve the MAX-CUT problem. QAOA takes place in the larger class of Variational Quantum Algorithms (VQAs)~\citep{Cerezo, Grange_VQA} which consist of alternating between a quantum circuit and a classical optimizer. Even if these algorithms have no performance guarantee for general problems, they are of great interest today because they have the convenient property of an adjustable quantum circuits' depth, making them implementable on the current NISQ computers~\citep{Preskill}. 

Recently, several combinatorial problems have been solved using QAOA. One can mention \emph{theoretical} problems such as MAX-CUT~\citep{QAOAFarhi}, Travelling Salesman Problem~\citep{TSP_Ruan}, MAX-3-SAT~\citep{Nannicini_Hyb_Perf}, Graph Coloring~\citep{GraphColoring}, Job Shop Scheduling~\citep{JobShop} and Vehicle Routing~\citep{Azad2022}. More \emph{industrial} problems have also been tackled, such as knapsack problem for battery revenue~\citep{Knapsack_Battery} or smart charging of electric vehicles~\citep{SmartCharging, Kea2023}. However, due to the small size of instances processed today (imposed by the weak maturity of quantum computers) and to the nature of heuristics whose performances are evaluated empirically, no quantum advantage is emerging yet.

Several levers are available when solving a problem with QAOA. For instance, one can choose the function that guides the classical optimizer, the depth or the type of gates used in the quantum circuit. Their impact in QAOA's performances have been proved~\citep{CVaR, Gibbs, Nannicini_Hyb_Perf}. 
However, the lever of the formulation of the combinatorial optimization problem has been less tackled in the literature, even though the shape of the quantum circuit is known to depend directly on the input problem. 
VQAs, and thus QAOA, necessitate the combinatorial optimization problem to be formulated as an Unconstrained Binary Optimization one, where the function is polynomial. While many real-world problems are constrained, transforming them into unconstrained problems represents a significant challenge, which could improve the performance of QAOA. The impact of the reformulation was first studied for Quadratic Unconstrained Binary Optimization (QUBO) reformulations. Indeed, QAOA derives from the Quantum Annealing algorithm of~\cite{Kadowaki98_QuantumAnnealing}, algorithm that takes as input an Ising Model, a model of ferromagnetism in statistical mechanics, which is equivalent to a QUBO problem. Thus, natively, the comparison of different formulations of QUBO models for a given problem recently raised interest in the community, for instance for the graph coloring problem~\citep{GraphColoring} or for Max-k-colorable subgraph problem~\citep{Quintero2022}. However, there are in fact no inherent degree limitations of the objective function. 
Some studies investigate Polynomial Unconstrained Binary Optimization (PUBO) formulations but, as far as we know, a generic reformulation approach is still lacking. These studies compare QAOA's performances of PUBO and QUBO formulations for specific problems such as for the Traveling Salesman Problem~\citep{Salehi2022}, the graph coloring problem~\citep{Campbell2022}, or a minimization problem with continuous variables~\citep{Stein2023}. However, they show an advantage of PUBO formulations that is often attributed to the reduction of the number of qubits required, though this may not be the only factor at stake. PUBO formulations offer several benefits: they not only reduce the number of qubits significantly but they also provide greater flexibility in shaping constraints, for instance by assigning specific binary values to the penalty functions. 

In this paper, we propose a generic reformulation for any polynomial problem into a PUBO problem, alongside with a specific reformulation in the case of linear constraints, exploiting PUBO’s ability to enforce penalty terms to take 0 and 1 values. We also present two generic QUBO reformulations for comparison, also applied to polynomial problems. Finally, we numerically compare QAOA performances of the PUBO formulation with binary penalty values, with one of the QUBO formulations on an industrial use-case: a railway timetabling problem. The timetabling we study consists of finding the transportation plan maximizing the operating profit according to the customers' demand taking into account the availability and cost of both the network and the rolling stock. It has been developed by the research department of the French Railway Company SNCF in order to help operators conceiving optimal high-speed train timetables. Its simplification, necessary because of the weak maturity of current quantum computers, leads to a Bin Packing problem with slight modifications on the constraints.

The paper is structured as follows. 
In Section~\ref{sec:generic_reformulation}, we propose generic reformulations of any problem with integer variables, polynomial objective function and polynomial constraints into PUBO and QUBO problems. We also provide a specific PUBO reformulation for the case of linear constraints, enforcing penalty terms to take binary values. In Section~\ref{sec:railway_problem}, we describe the nominal railway timetabling problem, we present a simplified version and we reformulate it as a PUBO and a QUBO following the methods of Section~\ref{sec:generic_reformulation}. Eventually, in Section~\ref{sec:num_results}, we illustrate the impact of the reformulation on QAOA's performances by solving several instances, revealing that PUBO with binary penalty values performs better than QUBO for our problem. 

\section{Generic reformulations}
\label{sec:generic_reformulation}
Variational Quantum Algorithms (VQAs) tackle unconstrained problems. However, most real-world combinatorial problems are constrained. Some constraints are directly related to the definition of the problem, for instance, that a city is visited exactly once by the salesman in the TSP. Some others express real-world limits, such as the limited number of seats on a train or the finite size of a knapsack, involving numerical constants in the description of the problem. In both cases, we need to reformulate the problem as an unconstrained optimization problem to solve it with VQAs. To remove the constraints, a common approach is to integrate them as penalty terms in the objective function. Examples of such reformulations for several NP-hard problems into QUBO problems are proposed in the literature~\citep{Glover,Lucas}. Notice that other techniques to remove constraints have been proposed, such as modifying the circuit of VQAs to \emph{express} the constraints~\citep{Hadfield}. 

In this section, we aim to widen the range of problems we can address with VQAs. For that, we provide a general method to transform any problem with integer variables, a polynomial objective function, and polynomial constraints into an \emph{equivalent} PUBO problem. Then, taking advantage of the wide possibilities of PUBO formulations, we propose a specific PUBO formulation with binary penalty terms, for problems with linear constraints only. Finally, we provide two different general transformations into an \emph{equivalent} QUBO problem. Two problems are said to be equivalent if and only if they have the same optimal solution(s) with the same optimal value.

\subsection{Integration of constraints}
First, we present a generic way to remove constraints from a nominal constrained problem and ensure that the reformulation into the resulting unconstrained problem is valid, i.e. that the resulting problem is \emph{equivalent} to the nominal problem. Specifically, let us consider the problem 
\begin{align}
\tag{\text{Nominal}}
\label{prob:initial_problem}
    & \begin{array}{rl}
    \min\limits_{x \in \mathcal{X}} & f(x) \\
    \text{subject to} & \text{the set of constraints } \{k \in \mathcal{K}\}
    \end{array}
\end{align}
For each constraint $k\in\mathcal{K}$, we define a penalty function $\pi_k : \mathcal{X} \rightarrow \R^+$ that satisfies, for $x\in \mathcal{X}$,
\begin{equation*}
    \pi_k(x) \left\{\begin{aligned}
      & = 0 &\text{if $x$ satisfies constraint $k$} \\
      & \geq 1 &\text{if $x$ violates constraint $k$}
    \end{aligned}\right.
\end{equation*}
We reformulate the constrained problem~\eqref{prob:initial_problem} as follows.
\begin{definition}
The unconstrained problem, for which we integrate the constraints of the nominal problem as penalty terms in the objective function, is
\begin{equation}
\tag{\text{Unconstrained}}
\label{prob:initial_problem_into_unconstrained_problem}
    \min_{x\in \mathcal{X}} f(x) + \sum_{k\in\mathcal{K}} \lambda_k \pi_k(x)\,,
\end{equation}
where the $\lambda_k > 0$ are penalty coefficients. 
\end{definition}
Thus, the reformulation of~\eqref{prob:initial_problem} into~\eqref{prob:initial_problem_into_unconstrained_problem} not only requires finding the penalty functions $\pi_k$, as we will discuss in Subsections~\ref{subsec:generic_PUBO_transformation} and~\ref{subsec:generic_QUBO_transformation}, but also requires choosing the numerical values of the penalty coefficients $\lambda_k$. Next, we provide a general lower bound for them.  

Let us note $f_\text{min} := \min\{f(x) : x\in\mathcal{X}\}$ the minimum value of $f$, $f_\text{max}:= \max \{f(x) : x\in\mathcal{X}\}$ its maximum value, and $f^* := \min\limits_{x\in\mathcal{X}} \{f(x) : x \mbox{ respects constraint } k,\forall k\in \mathcal{K}\}$ the optimal value of~\eqref{prob:initial_problem}. 
\begin{proposition}
    If we set, for all $k\in \mathcal{K}$, $$\lambda_k \geq f_\text{max} - f_\text{min}\,,$$ thus, we ensure that~\eqref{prob:initial_problem} is equivalent to~\eqref{prob:initial_problem_into_unconstrained_problem}.
\end{proposition}

\proof
The smallest value of an unfeasible solution of~\eqref{prob:initial_problem_into_unconstrained_problem} is always larger than the minimum value of $f$ plus the penalty cost of violating at least one constraint $k'$, namely, larger than
\begin{equation*}
    f_\text{min} + \lambda_{k'}\,.
\end{equation*}
Thus, for $\lambda_k \geq f_\text{max} - f_\text{min}$ for any constraint $k$, the smallest value of an unfeasible solution of~\eqref{prob:initial_problem_into_unconstrained_problem} is larger than
$$f_\text{min} + \lambda_{k'} \geq f_\text{min} + f_\text{max} - f_\text{min} = f_\text{max} > f^*\,.$$
Moreover, by definition of the penalty function, the values of the objective function of the two problems coincide on feasible solutions. Consequently, the optimal value of~\eqref{prob:initial_problem}, and its corresponding solution, is equal to the one of~\eqref{prob:initial_problem_into_unconstrained_problem}.
\endproof

Notice that if we set $\lambda_k \geq f_\text{max} - f_\text{min}$ for all $k\in \mathcal{K}$, not only~\eqref{prob:initial_problem_into_unconstrained_problem} has the same optimal solution as~\eqref{prob:initial_problem}, but also any feasible solution (of~\eqref{prob:initial_problem}) has a lower loss function value in~\eqref{prob:initial_problem_into_unconstrained_problem} than any unfeasible solution (of~\eqref{prob:initial_problem}). This latter property is worth noting because VQAs are heuristics, so they do not always find the optimal solution but solutions close to the optimal (in terms of the loss function). In practice, if $f_\text{min}\geq 0$, an upper bound of $f_\text{max}$ provides a lower bound for each $\lambda_k$. Thus, setting $\lambda_k \geq f_\text{max}$ ensures the new unconstrained problem to satisfy the above-mentioned properties.

In this subsection, we showed that reformulating a constrained problem into an unconstrained problem amounts to finding the penalty functions for each constraint. Next, we present a broad class of problems for which we provide the expression of the penalty functions.

\subsection{Class of eligible problems}
\label{subsec:eligible_class}
We present a class of problems~\eqref{prob:class_eligible_problems} for which we provide next a method to reformulate them as PUBO and QUBO problems. This class contains problems with integer variables, a polynomial objective function, and polynomial constraints. 

\begin{definition}
    Let $n\in \N$ be the number of variables and let $m\in \N$ be the number of constraints. We call~\eqref{prob:class_eligible_problems} the following class of problems.

\begin{align}
\label{prob:class_eligible_problems}
\tag{\text{IP-poly}}
& \min_{x}
    & & f(x_1,\ldots, x_n) \\
    & \textup{subject to} 
    & & g_k(x_1,\ldots, x_n) \leq 0\,, & \forall k \in [m] \nonumber\\
    &&& x_i \in \N, & \forall i \in [n] \nonumber
  \end{align}
where the functions $f$ and $g_k$, for any $k\in [m] :=\{1,\ldots, m\}$, are polynomial. Specifically,
\begin{align}
\tag{\text{IP-poly}}& \min_{x}
    & & \sum_{\gamma=(\gamma_1,\ldots, \gamma_n) \in \Gamma}{\alpha_{\gamma}x_1^{\gamma_1}\ldots x_n^{\gamma_n}} \\
    & \textup{subject to} 
    & & \sum_{\gamma=(\gamma_1,\ldots, \gamma_n) \in \Gamma_k}{\beta_{k,\gamma}x_1^{\gamma_1}\ldots x_n^{\gamma_n}} \leq 0\,,  & \forall k \in [m] \tag{$\text{C}_k$}\label{lign:Eligible_pb_constraint_k} \\
    &&& x_i \in \N, & \forall i \in [n] \nonumber
\end{align}
where $\Gamma \subseteq \N^n$ is a finite set and $\alpha_\gamma \in \R$ for $\gamma\in \Gamma$; and for all $k\in [m]$, $\Gamma_k\subseteq \N^n$ is a finite set and $\beta_{k,\gamma}\in \Z$ for $\gamma\in \Gamma_k$.
\end{definition}

\subsection{Reformulation into PUBO}
\label{subsec:generic_PUBO_transformation}
In this subsection, we present a generic method to reformulate any problem of~\eqref{prob:class_eligible_problems} into a PUBO problem (Subsubsection~\ref{subsubsec:generic_PUBO_transformation_general}). We also provide a reformulation in the case of linear constraints of the nominal problem (Subsubsection~\ref{subsubsec:generic_PUBO_transformation_binary_penalty_terms}), enforcing a certain shape of the penalty functions. 

\subsubsection{General PUBO reformulation}
\label{subsubsec:generic_PUBO_transformation_general}
The two steps to reformulate any problem of~\eqref{prob:class_eligible_problems} into a PUBO problem are the following. First, we transform the integer variables into binary variables. Second, we integrate the constraints into the objective function as penalty terms. Let us specify each of these steps.

The first step is the transformation of integer variables into binary variables. For that, we replace each integer variable $x_i$, for $i\in [n]$, by its binary decomposition $$x_i = \sum_{j=0}^{\lfloor \log_2(x) \rfloor}{x_i^{(j)}2^j}\,.$$
Thus, this decomposition requires $\lfloor \log_2(x) \rfloor + 1$ binary variables $x_i^{(j)} \in \{0,1\}$.

The second step is the integration of the constraints into the objective function as penalty terms. It is done as follows. 
Let us consider the constraint~\eqref{lign:Eligible_pb_constraint_k}
\begin{equation}
\tag{$\text{C}_k$}
    \sum_{\gamma=(\gamma_1,...\gamma_n) \in \Gamma_k}{\beta_{k,\gamma}x_1^{\gamma_1}...x_n^{\gamma_n}} \leq 0\,.
\end{equation}
This step aims at finding a penalty function for the above-mentioned constraint. 
Notice that after the first step, all variables are binary, so the integer variables of the left-hand side of~\eqref{lign:Eligible_pb_constraint_k} would be replaced by their binary description, adding more terms to the sum. Thus, without loss of generality, we assume henceforth that all the $x_i$ are binary variables. It results the following upper bound: $$\absolute{\sum_{\gamma \in \Gamma_k}{\beta_{k,\gamma}x_1^{\gamma_1}...x_n^{\gamma_n}}} \leq \sum_{\gamma\in \Gamma_k} \absolute{\beta_{k,\gamma}} =: \textup{UB}_k\,.$$ Thus, we define the penalty function associated with Constraint~\eqref{lign:Eligible_pb_constraint_k} as follows. 

\begin{proposition}
For $k\in [m]$, the function
    \begin{equation*}
    \pi_k^{\text{PUBO}}(x_1,\ldots,x_n) = \prod_{j=0}^{\textup{UB}_k} \left( \sum_{\gamma \in \Gamma_k}{\beta_{k,\gamma}x_1^{\gamma_1}...x_n^{\gamma_n}} + j \right)
\end{equation*}
is a penalty function for Constraint~\eqref{lign:Eligible_pb_constraint_k}. 
\end{proposition}

\proof
On the one hand, if $(x_1,\ldots,x_n)$ satisfies the constraint, it means that $\sum_{\gamma \in \Gamma_k}{\beta_{k,\gamma}x_1^{\gamma_1}...x_n^{\gamma_n}}$ takes a value in $\llbracket - \textup{UB}_k, 0\rrbracket$ and then there exists $j \in \llbracket 0, \textup{UB}_k\rrbracket$ that makes the product equal to 0, i.e. $\pi_k^{\text{PUBO}}(x_1,\ldots,x_n) = 0$. On the other hand, if $(x_1,\ldots,x_n)$ violates the constraint, the term $\sum_{\gamma \in \Gamma_k}{\beta_{k,\gamma}x_1^{\gamma_1}...x_n^{\gamma_n}}$ is strictly positive. Precisely, because each $\beta_{k,\gamma}$ is in $\Z$, the term $\sum_{\gamma \in \Gamma_k}{\beta_{k,\gamma}x_1^{\gamma_1}...x_n^{\gamma_n}}$ cannot be smaller than 1, leading to $\pi_k^{\text{PUBO}}(x_1,\ldots,x_n) \geq 1$.
\endproof

\subsubsection{Case of linear constraints}
\label{subsubsec:generic_PUBO_transformation_binary_penalty_terms}
Let us suppose that the constraints of~\eqref{prob:class_eligible_problems} are linear. In that case, we propose penalty functions that take binary values. Specifically, the penalty function takes value 1 if the corresponding constraint is violated and 0 otherwise. The choice of such binary functions is motivated by understanding wether of not controlling the violation of cost of each constraint improves performances of VQAs to solve such problems. This is discussed in Section~\ref{sec:num_results}.

In what follows, we provide penalty functions of binary values of equality and inequality constraints where the variable-dependent term is the sum of binary variables. Notice that for the case of a \emph{weighted} sum of binary variables (with weights in $\Z$), the penalty function can be easily deduced.

Let us begin with the penalty term for equality constraints. 

\begin{restatable}{property}{GenericCstEquality}
\label{prop:generic_cst_equality}
Let us consider the constraint, for $n\in\N$ and $c\in\N$,
\begin{equation}
\label{eq:constraint_eq}
\tag{Eq}
    \sum_{i=1}^n x_i = c\,.
\end{equation}
For $c\in\N^*$, the penalty function is, for $x\in\{0,1\}^n$,
\begin{equation*}
    \piEq_c(x) = 1 + \sum_{k=c}^n (-1)^{k-c+1} \binom{k}{c}\sum_{i = \{i_1,\ldots, i_k\} \in \EnsBinomial{n}{k}} x_{i_1}\ldots x_{i_k}\,,
\end{equation*}
where $ \EnsBinomial{n}{k}$ denotes all the sets of $k$ elements in $[n]$.
For the specific case of $c=0$, the penalty function is
\begin{equation*}
    \piEq_0(x) = \sum_{k=1}^n (-1)^{k+1}\sum_{i \in \EnsBinomial{n}{k}} x_{i_1}\ldots x_{i_k}\,.
\end{equation*}
\end{restatable}

See proof in Appendix~\ref{appendix:proof_cst_equality}. Next, we define a penalty term for inequality constraints. Property~\ref{prop:generic_cst_inequality_inf} deals with the inferiority case whereas Property~\ref{prop:generic_cst_inequality_sup} tackles the superiority case.

\begin{restatable}{property}{GenericCstInequalityInf}
\label{prop:generic_cst_inequality_inf}
We consider the constraint, for $n\in\N$ and $c\in\N$,
\begin{equation}
\label{eq:constraint_inf}
\tag{Inf}
    \sum_{i=1}^n x_i \leq c\,.
\end{equation}
The associated penalty function is, for $x\in\{0,1\}^n$,
\begin{equation*}
    \piInf_c(x) = \sum_{k=c+1}^n (-1)^{k-c+1} \binom{k-1}{c}\sum_{i\in \EnsBinomial{n}{k}} x_{i_1}\ldots x_{i_k}\,.
\end{equation*}
\end{restatable}

See proof in Appendix~\ref{appendix:proof_cst_inequality_inf}. 

\begin{restatable}{property}{GenericCstInequalitySup}
\label{prop:generic_cst_inequality_sup}
We consider the constraint, for $n\in\N$ and $c\in\N^*$,
\begin{equation}
\label{eq:constraint_sup}
\tag{Sup}
    \sum_{i=1}^n x_i \geq c\,.
\end{equation}
The associated penalty function is, for $x\in\{0,1\}^n$,
\begin{equation*}
    \piSup_c(x) = 1 + \sum_{k=c}^n (-1)^{k-c+1} \binom{k-1}{c-1}\sum_{i\in \EnsBinomial{n}{k}} x_{i_1}\ldots x_{i_k}\,.
\end{equation*}
\end{restatable}

See proof in Appendix~\ref{appendix:proof_cst_inequality_sup}.

\subsection{Reformulation into QUBO}
\label{subsec:generic_QUBO_transformation}
Any problem of~\eqref{prob:class_eligible_problems} can also be reformulated as a QUBO problem. We present two possible generic ways to do so.

The first one is to compute the PUBO penalty term $\pi_k^{\text{PUBO}}$ resulting of the two steps presented in the previous subsection, and hereafter to decrease its degree until it is equal to 2. For that, we reduce recursively the degree of each monome of degree larger than 3 by using the following method of linearization.
Repeat until no monome of degree at least 3 exists: 
\begin{enumerate}
    \item Choose a monome of degree larger than 3.
    \item Pick two (binary) variables appearing in the monome $x_i$ and $x_j$ and replace the product $x_ix_j$ by the new binary variable $y_{ij}$ in every monomes of the penalty function $\pi_k^{\text{PUBO}}(x)$ where it appears.
    \item Add the term $\lambda(x_ix_j - 2x_iy_{ij} - 2 x_jy_{ij} + 3y_{ij})$ to $\pi_k^{\text{PUBO}}$, with $\lambda>0$ a constant.
\end{enumerate}
Thus, it decreases the degree of the considered monome (and possibly other monomes) by 1. Notice that the term $$pen_\text{lin}(x_i,x_j,y_{ij}) := x_ix_j - 2x_iy_{ij} - 2 x_jy_{ij} + 3y_{ij}$$ is a penalty term associated to the constraint linearization $x_ix_j = y_{ij}$. Indeed, Table~\ref{tab:truth_table} shows that $pen_\text{lin}$ is equal to zero if the value of the triplet $(x_i,x_j,y_{ij})$ satisfies the constraint (rows in bold), and larger than 1 otherwise.

\begin{table}
    \centering
    \begin{tabular}{cccc}
        $x_i$ & $x_j$ & $y_{ij}$ & $pen_\text{lin}(x_i,x_j,y_{ij})$ \\
        \hline
         $\textbf{0}$ & $\textbf{0}$ & $\textbf{0}$ & $\textbf{0}$\\
         $\textbf{0}$ & $\textbf{1}$ & $\textbf{0}$ & $\textbf{0}$\\
         $\textbf{1}$& $\textbf{0}$ & $\textbf{0}$ & $\textbf{0}$\\
         1& 1 & 0 & 1\\
         0& 0 & 1 & 3\\
         0& 1 & 1 & 1\\
         1& 0 & 1 & 1\\
         $\textbf{1}$& $\textbf{1}$ & $\textbf{1}$ & $\textbf{0}$\\
    \end{tabular}
    \caption{Values of function $pen_\text{lin}$}
    \label{tab:truth_table}
\end{table}

The second way to provide quadratic penalty terms is the following. Given the constraint with binary variables~\eqref{lign:Eligible_pb_constraint_k} resulting from the first step of PUBO formulation presented in the previous subsection, linearize it as proposed above. Let us note $\text{Lin}_k(x,y) \leq 0$ the resulting constraint, where the $y$ variables are resulting from the linearization process. Thus, the penalty term associated with this constraint is $$\pi_k^{\text{QUBO}}(x,y) = \min_{s\in\llbracket0,-\min_{x,y} \textup{Lin}_k(x,y) \rrbracket}(\textup{Lin}_k(x,y) + s)^2\,.$$ Indeed, for a given $(x,y)$:
\begin{itemize}
    \item If the constraint is satisfied, i.e. $\text{Lin}_k(x,y) \leq 0$, there exists a value of $s \in\llbracket0,-\min\limits_{x,y} \textup{Lin}_k(x,y) \rrbracket$ such that $\text{Lin}_k(x,y) + s = 0$, and thus $\pi_k^{\text{QUBO}}(x,y) = 0$.
    \item Otherwise, $\text{Lin}_k(x,y) > 0$ and there is no value in $s \in\llbracket0,-\min_{x,y} \textup{Lin}_k(x,y) \rrbracket$ such that $s = - \textup{Lin}_k(x,y)$, and thus $\pi_k^{\text{QUBO}}(x,y) \geq 1$.
\end{itemize}
Notice that we suppose that the constraint can be satisfied, namely that $\min_{x,y} \textup{Lin}_k(x,y)\leq 0$.

Eventually, because VQAs can only deal with one optimization problem, the optimization problems $\min_s$ resulting from the expression of the penalty functions $\pi_k^{\text{QUBO}}$ must join the initial optimization problem $\min_{x,y}$ of the loss function resulting of the first step of PUBO formulation. Gathering the optimization problems of penalty terms with the nominal optimization problem leads to a overall minimization problem over the slack variables $(s)$ in addition to the decision variables $(x,y)$. The optimal solution(s) of the overall minimization problem coincide with the optimal solution(s) of the two-level optimization problem, but it is worth noting that the optimization process is different and leads, among other, to non-optimal values of slack variables for optimal values of decision variables. This phenomenon is discussed in Section~\ref{sec:num_results}.

Thus, there are $m$ additional variables $s_k \in\llbracket0,-\min_{x,y} \textup{Lin}(x,y) \rrbracket$ (eventually written with binary variables), for $k\in[m]$, in the overall optimization QUBO problem. Henceforth, we call these variables \emph{slack} variables to emphasize the fact that they do not represent nominal decision variables but additional variables coming from the reformulation. 

Next, we present the railway timetabling problem at hand and apply both QUBO and PUBO reformulations presented in the section.

\section{A railway timetabling problem and its reformulations}
\label{sec:railway_problem}

\subsection{Nominal problem}
Railway timetabling problems are crucial problems for railway companies. A first version of the timetable is often planned several years in advance, and is related to other planning problems, such as crew scheduling or rolling stock scheduling. The goal of timetabling is to ensure the satisfaction of customers, the minimization of delays thanks to robustness and resilience properties, the minimization of costs, and the validation of both operational and security constraints. 

The railway timetabling problem for high-speed trains at the French Railway Company SNCF we are considering is the following: according to the customers' demand (estimated from past data) and the availability of the network and the rolling stock, the aim is to find the transportation plan maximizing the operating profit. The output is a timetable of trains, the associated rolling stock schedule, and a forecast of the passengers for each train and journey. The optimal solution is the best compromise between the revenues generated by the customers' journeys and the production costs (network, rolling stock, human resources, etc.). It is formulated as an Integer Linear Programming (ILP). 
Let us provide a basic description of it while explaining the important ideas to understand the simplification proposed below.
First, we describe the main sets of an instance of our problem. 
\begin{itemize}
    \item S, the set of Train-paths. A train-path is a timed unitary portion of tracks. It defines the availability to run a carriage over a portion of tracks over a given time period.
    \item T, the set of available Trains. A train is described as a union of train-paths. A train is defined by its origin, destination, and the served stations, with departure and arrival times for each station. A train uses the same carriage throughout the journey.
    \item G, the set of Groups of customers. A group gathers customers that have the same preferences on journeys, namely, customers wanting to leave, respectively arrive, at the same station and at the same time.
    \item $R =\{(t,g)\in T\cross G : \text{group $g$ accepts to take train $t$}\}$, the set of possible Customers' journeys. This set expresses the possibilities to satisfy the customer groups' demands.
\end{itemize}
Other sets are required, such as the set expressing the different incompatibilities between train-paths or the set of trains that can be coupled. Besides sets, many constants appear in the original formulation, both in the objective function and in the constraints, such as the maximum capacity of a carriage, the toll cost of a train-path, or the average receipt for a journey. 

Second, let us introduce the variables of this problem, which are binary or integers.
\begin{itemize}
    \item $x\in\{0,1\}^{|T|}$, where $x_t = 1$ iff $t\in T$ is used in the timetable.
    \item $u\in\{0,1\}^{|S|}$, where $u_s = 1$ iff the train-path $s\in S$ is used.
    \item $z\in \mathbb{N}^{|R|}$, where $z_r$ is the number of customers for journey $r\in R$.
\end{itemize}

The objective function leads to finding the timetable providing the best compromise between customers’ demand satisfaction and production cost. Specifically,
\begin{equation*}
    \max_{x_t, u_s, z_r} \alpha(z_r) - \beta(x_t,u_s)\,,
\end{equation*}
where $\alpha(z_r)$ is a linear function representing the receipts generated by selling tickets and 
$\beta(x_t,u_s)$ is a linear function representing the total cost associated with the use of train-paths and carriages. 
An optimal solution to our problem is a feasible timetable that maximizes the above loss function. The feasibility of a timetable is defined by many linear constraints such as forbidding customers to take a train not used in the timetable, ensuring that the maximum capacity of a train is satisfied on each train-path, or setting a minimum daily frequency for a given journey.  
In practice, this railway timetabling problem requires about twenty sets and constants to define an instance of the nominal problem, while it requires four types of binary or integer variables, a linear function, and ten types of linear constraints, two of which are soft constraints. Today, this problem is solved with an ILP solver and, depending on the instance, it can be hard or even impossible to find an optimal or near-optimal solution. For example, the optimal solution of the problem for the sector between two French cities Paris and Lyon is found within a second whereas the solution for timetabling of the inter-regional trains has a gap of 67\% from the optimal solution after 10 minutes of running time. No solution is found for larger instances such as the perimeter of the entire metropolitan French territory. 

While Section~\ref{sec:generic_reformulation} proves that this railway timetabling problem can be solved theoretically with quantum-classical metaheuristics, because it belongs to~\eqref{prob:class_eligible_problems}, we need to simplify it by putting aside some assumptions to be able to solve it on current quantum machines, or classical simulators of quantum machines. Indeed, the size of the original problem would be too large even considering instances on small geographical sectors and short periods of time (one day). For example, the instance on the sector Paris-Lyon for one day only, which covers 6 stations, amounts to 176 customer groups for around 25 000 customers, 87 feasible trains, and 157 train-paths leads to a nominal problem with 8 000 binary variables. Additionally, VQAs require unconstrained problems (QUBO or PUBO problems), integrating constraints in the objective function costs in terms of additional variables (for QUBO reformulation) and additional gates (for QUBO and PUBO reformulation) as explained in Section~\ref{sec:generic_reformulation}. As an example, the QUBO formulation of the above-mentioned instance requires roughly 12 000 additional qubits, ending up to 20 000 qubits for the total description of the instance. This size prohibits a resolution on current gate-based quantum hardware that does not exceed a hundred qubits with a reasonable gate error rate, without mentioning the high connectivity that would be necessary. While keeping the essence of the initial timetabling problem, we consider a simplification that is an extended version of a Bin Packing problem, as we detail next.

\subsection{Extended Bin Packing Problem simplification}


Let us present the simplified version of the nominal railway problem that results in an Extended Bin Packing problem. This simplification retains the core aspects of the original problem while relaxing certain operational requirements. 
As mentioned above, this simplification is motivated by two things. The first one is to reduce the number of decision variables, and thus the number of qubits, to describe an instance. The second one is to ease the transformation of the constrained problem into an unconstrained one, avoiding too many additional qubits and too many gates for the implementation of the quantum circuit of VQAs. Moreover, the Extended Bin Packing problem is NP-hard as the nominal problem, which comforts the interest of this choice. 

The simplification is done as follows. Let $m\in\N$ be the number of customer groups and $n\in\N$ the number of available trains. We define
\begin{itemize}
    \item $G := \{g_1,\ldots,g_m\}$, a set of $m\in\N$ customer groups.
    \item $T := \{T_1,\ldots,T_n\}$, the set of $n\in\N$ available trains, where each train $T_i\subseteq G$ is the set of groups that accept to take this train.
\end{itemize}
We consider that a customer group $g_k$, for $k\in [m]$, is satisfied if at least one of the trains matching its demand (the $T_i$s such that $g_k\in T_i$) is in the output timetable. 

Additionally, we suppose that each group contains the same number of customers. The latter assumption does not cause too much loss of generality because one can always define the smallest group as a unit and duplicate groups that are bigger. For each available train $T_i$, $i\in[n]$, we specify
\begin{itemize}
    \item $p_i\in\R^+$, the benefit of selling tickets to one group for train $T_i$.
    \item $c_i\in\R^+$, the cost of using train $T_i$.
\end{itemize}
We note $\CMax\in\N$ the maximum number of groups a carriage can accommodate, i.e. the maximum capacity of a carriage divided by the (fixed) number of customers in a group. We define below the problem considered, called Extended Bin Packing.

\begin{definition}[Extended Bin Packing problem]
The binary decision variables for the Extended Bin Packing problem are of two types. The first one indicates if the train is taken in the timetable: $\forall i\in [n]$,
\begin{equation*}
    x_i = \begin{cases}
    1 \hspace{2mm} \textup{if $T_i$ is taken in the timetable}\\
    0 \hspace{2mm} \textup{else}
    \end{cases}
\end{equation*}
The second assigns groups to trains in the timetable: $\forall i\in [n], \forall j \in [m]$,
\begin{equation*}
    y_{ij} = \begin{cases}
    1 \hspace{2mm} \textup{if group $g_j$ takes train $T_i$}\\
    0 \hspace{2mm} \textup{else}
    \end{cases}
\end{equation*}
Notice that we declare a variable $y_{ij}$ if and only if $g_j\in T_i$, namely that train $T_i$ can satisfy group $g_j$. It avoids unnecessary variables and reduces the number of $y$ variables from $nm$ to $q := \sum_{i=1}^n |T_i|$. The problem can be seen as a variation of the Bin Packing problem, where each bin has a different capacity, and not every item has to be put in a bin because we are rather looking for the best compromise between the cost of using bins and the reward of putting items into bins. Hence we call this modelization the Extended Bin Packing problem. It is stated as follows: 

\begin{align}
\tag{Extended-BP}
\label{problem:SCP_extended}
    &\min_{x,y}
    & &\sum_{i=1}^n c_ix_i - \sum_{i=1}^n\sum_{\substack{j\in[m] :\\g_j\in T_i}} p_iy_{ij}  \\
    & \textup{subject to} 
    & & \sum_{\substack{i\in[n]:\\g_j\in T_i}} y_{ij} \leq 1, & \forall j \in [m]\tag{\text{Uni}}\label{lign:satisfaction}\\
    &&& \sum_{\substack{j\in[m]:\\g_j\in T_i}} y_{ij} \leq \CMax\cdot x_i, & \forall i \in [n]\tag{\text{Capa}}\label{lign:capacity}\\
    &&& x_i \in \{0,1\}, & \forall i \in [n] \nonumber\\
    &&&  y_{ij} \in \{0,1\}, & \forall i \in [n],\, j \in [m] \text{ such that } g_j\in T_i\nonumber
  \end{align}
Constraint~\eqref{lign:satisfaction} ensures that each customer group takes at most one train from those available in the timetable, and Constraint~\eqref{lign:capacity} both expresses the limited capacity of a carriage and forbids a group to take a train unused in the timetable.
\end{definition}

Next, we provide reformulations of this simplified model into QUBO and PUBO problems to be able to solve it with QAOA. 

\subsection{Reformulations into QUBO and PUBO problems}

In this subsection, we provide reformulations of~\eqref{problem:SCP_extended} into PUBO and QUBO problems. For the PUBO formulation, we apply the method presented in Subsubsection~\ref{subsubsec:generic_PUBO_transformation_binary_penalty_terms} which has the specificity to have binary penalty values. For the QUBO formulation, we apply the second method of Subsection~\ref{subsec:generic_QUBO_transformation}, which is the most common QUBO reformulation encountered in the literature.
For both reformulation, we need to integrate two different types of constraints,~\eqref{lign:satisfaction} and~\eqref{lign:capacity}. We recall them below. The first constraint, expressing that a group takes at most one train, is: for all $j\in [m]$, 
\begin{equation}
\label{cst:Uni_j}
\tag{$\text{Uni}_j$}
     \sum_{\substack{i\in[n]:\\g_j\in T_i}} y_{ij} \leq 1\,.
\end{equation}

The second type of constraint, the capacity constraint, is: for all $i\in[n]$, 
\begin{equation}
\label{cst:Capa_i}
\tag{$\text{Capa}_i$}
     \sum_{\substack{j\in[m]:\\g_j\in T_i}} y_{ij} \leq \CMax\cdot x_i\,.
\end{equation}

\subsubsection{PUBO reformulation with penalty binary values}
 

Based on Subsection~\ref{subsec:generic_PUBO_transformation}, and more specifically on the properties of Subsubsection~\ref{subsubsec:generic_PUBO_transformation_binary_penalty_terms}, it results the following penalty terms for the Extended Bin Packing problem as a PUBO.

\begin{proposition}[Penalty term for Constraint~\eqref{cst:Uni_j}]
\label{prop:cst_Uni_j}
    Let $j\in [m]$ and let us note $G_j$ the set of indices of trains accepted by group $g_j$. The function, for $x\in\{0,1\}^n$, $y\in \{0,1\}^q$, \begin{equation*}
    penUni_j(x,y) = \sum_{k=2}^{|G_j|}(-1)^k(k-1)\sum_{\substack{I = (i_1,\ldots,i_k)\subseteq G_j :\\ |I|=k}}y_{i_1,j}\ldots y_{i_k,j}
\end{equation*} 
is a penalty term for Constraint~\eqref{cst:Uni_j}.
\end{proposition}

\proof
Use Property~\ref{prop:generic_cst_inequality_inf} for the binary variables $y_{ij}$ such that $i\in G_j$, and for the constant $c=1$. 
\endproof

\begin{proposition}[Penalty term for Constraint~\eqref{cst:Capa_i}]
\label{prop:cst_Capa_i}
    Let $i\in [n]$. The following function is a penalty term for Constraint~\eqref{cst:Capa_i}. For $x\in\{0,1\}^n$, $y\in \{0,1\}^q$,
    \begin{equation*}
    penCapa_i(x,y) = (1-x_i)\cdot pen_0(x,y) + x_i\cdot pen_1(x,y)\,,
\end{equation*}
where
\begin{equation*}
    pen_0(x,y) = \sum_{k=1}^{|\{j:i\in G_j\}|}(-1)^{k+1}\sum_{\substack{J\subseteq \{j:i\in G_j\}:\\ |J|=k}} y_{i,j_1}\ldots y_{i,j_k}
\end{equation*}
and 
\begin{equation*}
    pen_1(x,y) = \sum_{k=\CMax + 1}^{|\{j:i\in G_j\}|}(-1)^{k-\CMax +1}\binom{k-1}{\CMax}\sum_{\substack{J\subseteq \{j:i\in G_j\} :\\ |J|=k}} y_{i,j_1}\ldots y_{i,j_k}\,.
\end{equation*}
\end{proposition}

\proof
We distinguish two cases. 
\begin{itemize}
    \item If $x_i = 0$, then we use Property~\ref{prop:generic_cst_equality} for the binary variables $y_{ij}$ such that $j\in\{j:i\in G_j\}$, and for the constant $c=0$, which provides the term $pen_0(x,y)$. 
    \item If $x_i = 1$, then we use Property~\ref{prop:generic_cst_inequality_inf} for the binary variables $y_{ij}$ such that $j\in\{j:i\in G_j\}$, and for the constant $c=\CMax$, which provides the term $pen_1(x,y)$. 
\end{itemize}
Eventually, we add the two terms as follows, each one appearing when the condition on $x_i$ is satisfied, to obtain the final penalty function: $(1-x_i)\cdot pen_0(x,y) + x_i\cdot pen_1(x,y)$. 
\endproof

It results from the two propositions above the reformulation of our simplified problem~\eqref{problem:SCP_extended} into a PUBO problem. 

\begin{proposition}[Extended-BP into PUBO]
The reformulation of~\eqref{problem:SCP_extended} as a PUBO problem is:
\begin{align}
\label{problem:SCP_extended_PUBO}
\tag{Extended-BP-PUBO}
    \min\limits_{x\in\{0,1\}^n \atop y\in\{0,1\}^q} &\sum_{i=1}^n c_ix_i - \sum_{i=1}^n\sum_{\substack{j\in[m]:\\g_j\in T_i}} p_iy_{ij}\\ &+ \sum_{j=1}^m \lambda_{u,j} penUni_j(x,y) + \sum_{i=1}^n \lambda_{c,i}penCapa_i(x,y)\,,\nonumber
\end{align}
where $\lambda_{u,j}, \lambda_{c,i}\in \R^*_+$ for all $j\in[m], i\in[n]$ are the penalty coefficients, and where the penalty functions are defined in Propositions~\ref{prop:cst_Uni_j} and~\ref{prop:cst_Capa_i}.
\end{proposition}

\subsubsection{QUBO reformulation}
Applying the second method of Subsection~\ref{subsec:generic_QUBO_transformation}, the penalty terms for the Extended Bin Packing problem as a QUBO are the following.

\begin{proposition}[Quadratic penalty terms of Constraints~\eqref{cst:Uni_j} and~\eqref{cst:Capa_i}]
\label{prop:cst_uni_capa_qubo}
The penalty term for the first Constraint~\eqref{cst:Uni_j} is, for $j\in[m]$,
\begin{equation*}
    penUniQubo_j(x,y) = \min_{s_j\in\{0,1\}}\left( \sum_{i\in G_j}y_{ij} + s_j - 1  \right)^2\,.
\end{equation*}
The penalty term for the second constraint~\eqref{cst:Capa_i} is, for $i\in[n]$,
\begin{equation*}
    penCapaQubo_i(x,y) = \min_{r_i\in\llbracket0,\CMax\rrbracket}\left( \sum_{j:i\in G_j}y_{ij} - \CMax\cdot x_i + r_i  \right)^2\,.
\end{equation*}
\end{proposition}

 In the binary model, the integer slack variables $r_i$ are replaced by their binary expressions $r_i^{bin} := \sum_{l=0}^{\lfloor \log_2(\CMax)\rfloor} 2^k\cdot r_{i,k}^{bin}$, for $r_{i,k}^{bin} \in \{0,1\}$. In total, the penalty constraints introduce roughly $(m + n\cdot \log_2(\CMax))$ slack binary variables (denoted by $s$ and $r$) to the initial Extended Bin Packing problem to express the following QUBO formulation.

 \begin{proposition}[Extended-BP into QUBO]
 The reformulation of~\eqref{problem:SCP_extended} as a QUBO problem is
 \begin{align}
\label{problem:SCP_extended_QUBO}
\tag{Extended-BP-QUBO}
    \min_{\substack{x\in\{0,1\}^n,\\ y\in\{0,1\}^q \\ s\in \{0,1\}^m,\\ r^{bin}\in\{0,1\}^{n\cdot \log_2(\CMax)}}} &\sum_{i=1}^n c_ix_i - \sum_{i=1}^n\sum_{\substack{j\in[m]:\\g_j\in T_i}} p_iy_{ij}\\ &+ \sum_{j=1}^m \lambda_{u,j} \left( \sum_{i\in G_j}y_{ij} + s_j - 1  \right)^2 \\ &+ \sum_{i=1}^n \lambda_{c,i}\left( \sum_{j:i\in G_j}y_{ij} - \CMax\cdot x_i + r_i^{bin}  \right)^2\,, \nonumber
\end{align}
where $\lambda_{u,j}, \lambda_{c,i}\in \R^*_+$ for all $j\in[m], i\in[n]$ are the penalty coefficients.
\end{proposition}

Not only does the proposed QUBO formulation require additional qubits compared to the PUBO formulation, but it also has penalty terms that take integer values (and not binary values). These two facts are observed and discussed in the next section. 


\section{Numerical results}
\label{sec:num_results}

In this section, we illustrate the importance of the input problem formulation. Specifically, we compare the results of solving the~\eqref{problem:SCP_extended_QUBO} and the~\eqref{problem:SCP_extended_PUBO} reformulations with QAOA on several small instances, showing a trend regarding these two reformulations. For details on QAOA description and implementation, readers can refer to the paper of~\cite{Grange_VQA}. 

\subsection{Instances and results}
Let us consider the three instances below (Instance A, Instance B, Instance C). Each set represents a train and each element represents a customer group. In Table~\ref{tab:description_instances}, we provide for each instance: the numerical values for the problem description (cost of using a train $c$, benefit of selling a ticket $p$, maximum capacity $\CMax$), the optimal solutions and their respective optimal values, and the numerical values of the penalty coefficients (coefficients $\lambda_{u}, \lambda_{c}$) for the QUBO and PUBO reformulations. Notice that we only describe the set of trains in the optimal timetable, without specifying the assignment of each customer. The assignment is unique and straightforward for Instances A and B, but for Instance C, there are actually 11 optimal solutions, but still with the 3 sets of trains described in Table~\ref{tab:description_instances}. 

\begin{figure}[H]
\begin{minipage}{0.3\textwidth}
    \centering
    \resizebox{0.6\textwidth}{!}{\begin{tikzpicture}
    \draw[thick] (0,0) circle (1.25cm) (3,0) circle (1.25cm);
    \draw[fill=black] (0,0) circle (0.2cm) (3,0) circle (0.2cm);
    \draw[thick] (0,2) arc [radius=2, start angle=90, end angle=270] (3, -2) arc [radius=2, start angle=270, end angle=450] (0, 2) -- (3,2) (0,-2) -- (3,-2);
\end{tikzpicture}}
    \caption{Instance A}
    \label{fig:instanceA_SCPExtended}
\end{minipage}
\begin{minipage}{0.3\textwidth}
    \centering
    \resizebox{0.6\textwidth}{!}{\begin{tikzpicture}
    \draw[fill=black] (0,0) circle (0.2cm) (3,0) circle (0.2cm) (0,3) circle (0.2cm) (3,3) circle (0.2cm);
    \draw[thick] (0,1.25) arc [radius=1.25, start angle=90, end angle=270] (3,-1.25) arc [radius=1.25, start angle=-90, end angle=90] (0, 1.25) -- (3, 1.25) (0, -1.25)--(3,-1.25);
    \draw[thick] (0,4.25) arc [radius=1.25, start angle=90, end angle=270] (3,1.75) arc [radius=1.25, start angle=-90, end angle=90] (0, 4.25) -- (3, 4.25) (0, 1.75)--(3,1.75);
    \draw[thick] (-0.707, 0.707) arc [radius=1, start angle=135, end angle=315] (0.707, -0.707)--(3.707, 2.293) (3.707, 2.293) arc [radius=1, start angle=-45, end angle=135] (2.293, 3.707) -- (-0.707, 0.707);
\end{tikzpicture}}
    \caption{Instance B}
    \label{fig:instanceB_SCPExtended}
\end{minipage}
\begin{minipage}{0.3\textwidth}
    \centering
    \resizebox{0.6\textwidth}{!}{\begin{tikzpicture}
    \draw[fill=black] (0,0) circle (0.2cm) (-3, -2) circle (0.2cm) (3, 2) circle (0.2cm) (3, -2) circle (0.2cm) (-3, 2) circle (0.2cm);
    \draw[thick] (-2, 0 )arc [radius=3.5, start angle=180, end angle=540] ;
    \draw[thick] (2, 0 ) arc [radius=3.5, start angle=0, end angle=360];
    \draw[thick] (4.5,2) arc [radius=1.5, start angle=0, end angle = 180] (1.5,2) -- (1.5, -2) (1.5, -2) arc [radius=1.5, start angle=180, end angle =360] (4.5, -2) -- (4.5, 2);
\end{tikzpicture}}
    \caption{Instance C}
    \label{fig:instanceC_SCPExtended}
\end{minipage}
\end{figure}

\begin{table}[h!]
\small
    \centering
    \begin{tabular}{c|c|c|c} 
         & Instance A & Instance B & Instance C \\ \hline
        Number of trains & 3 & 3 & 3 \\ \hline
        Number of customer groups & 2 & 4 & 5 \\ \hline
        Cost $c_i, \forall i\in [n]$ & 1 & 1 & 1 \\ \hline
        Benefit $p_i, \forall i\in [n]$ & 1 & 1 & 1 \\\hline
        Maximum capacity $\CMax$ & 2 & 2 & 2 \\
        \hline \hline 
        Optimal solutions (in red) & \resizebox{14mm}{8mm}{\begin{tikzpicture}
    \draw[thick, black!50] (0,0) circle (1.25cm) (3,0) circle (1.25cm);
    \draw[fill=black] (0,0) circle (0.2cm) (3,0) circle (0.2cm);
    \draw[thick, red] (0,2) arc [radius=2, start angle=90, end angle=270] (3, -2) arc [radius=2, start angle=270, end angle=450] (0, 2) -- (3,2) (0,-2) -- (3,-2);
\end{tikzpicture}} & \resizebox{12mm}{12mm}{\begin{tikzpicture}
    \draw[fill=black] (0,0) circle (0.2cm) (3,0) circle (0.2cm) (0,3) circle (0.2cm) (3,3) circle (0.2cm);
    \draw[thick, red] (0,1.25) arc [radius=1.25, start angle=90, end angle=270] (3,-1.25) arc [radius=1.25, start angle=-90, end angle=90] (0, 1.25) -- (3, 1.25) (0, -1.25)--(3,-1.25);
    \draw[thick, red] (0,4.25) arc [radius=1.25, start angle=90, end angle=270] (3,1.75) arc [radius=1.25, start angle=-90, end angle=90] (0, 4.25) -- (3, 4.25) (0, 1.75)--(3,1.75);
    \draw[thick, black!50] (-0.707, 0.707) arc [radius=1, start angle=135, end angle=315] (0.707, -0.707)--(3.707, 2.293) (3.707, 2.293) arc [radius=1, start angle=-45, end angle=135] (2.293, 3.707) -- (-0.707, 0.707);
\end{tikzpicture}} & \makecell{\resizebox{15mm}{11mm}{\begin{tikzpicture}
    \draw[fill=black] (0,0) circle (0.2cm) (-3, -2) circle (0.2cm) (3, 2) circle (0.2cm) (3, -2) circle (0.2cm) (-3, 2) circle (0.2cm);
    \draw[thick, black!50] (-2, 0 )arc [radius=3.5, start angle=180, end angle=540] ;
    \draw[thick, red] (2, 0 ) arc [radius=3.5, start angle=0, end angle=360];
    \draw[thick, red] (4.5,2) arc [radius=1.5, start angle=0, end angle = 180] (1.5,2) -- (1.5, -2) (1.5, -2) arc [radius=1.5, start angle=180, end angle =360] (4.5, -2) -- (4.5, 2);
\end{tikzpicture}}\\ \resizebox{15mm}{11mm}{\begin{tikzpicture}
    \draw[fill=black] (0,0) circle (0.2cm) (-3, -2) circle (0.2cm) (3, 2) circle (0.2cm) (3, -2) circle (0.2cm) (-3, 2) circle (0.2cm);
    \draw[thick, red] (-2, 0 )arc [radius=3.5, start angle=180, end angle=540] ;
    \draw[thick, red] (2, 0 ) arc [radius=3.5, start angle=0, end angle=360];
    \draw[thick, red] (4.5,2) arc [radius=1.5, start angle=0, end angle = 180] (1.5,2) -- (1.5, -2) (1.5, -2) arc [radius=1.5, start angle=180, end angle =360] (4.5, -2) -- (4.5, 2);
\end{tikzpicture}} \\
        \resizebox{15mm}{11mm}{\begin{tikzpicture}
    \draw[fill=black] (0,0) circle (0.2cm) (-3, -2) circle (0.2cm) (3, 2) circle (0.2cm) (3, -2) circle (0.2cm) (-3, 2) circle (0.2cm);
    \draw[thick, red] (-2, 0 )arc [radius=3.5, start angle=180, end angle=540] ;
    \draw[thick, red] (2, 0 ) arc [radius=3.5, start angle=0, end angle=360];
    \draw[thick, black!50] (4.5,2) arc [radius=1.5, start angle=0, end angle = 180] (1.5,2) -- (1.5, -2) (1.5, -2) arc [radius=1.5, start angle=180, end angle =360] (4.5, -2) -- (4.5, 2);
\end{tikzpicture}}} \\ \hline
        Optimal value ($f^*$) & -1 & -2 & -2 \\
        \hline 
        Range of objective function values $(f_\text{min};f_\text{max})$ & $(-4;3)$ & $(-6;3)$ & $(-8;3)$ \\ \hline \hline 
        QUBO penalty coefficients $\lambda_{u,i} = \lambda_{c,j}, \forall i\in [n], j\in[m]$ & 8 & 10 & 12 \\ \hline
        PUBO penalty coefficients $\lambda_{u,i} = \lambda_{c,j}, \forall i\in [n], j\in[m]$ & 8 & 10 & 12 \\
        
    \end{tabular}
    \caption{Description of instances, their optimal solutions and the reformulation coefficients.}
    \label{tab:description_instances}
\end{table}

\begin{figure}[!h]
    \centering
    \subfigure{
        \includegraphics[width=0.46\textwidth]{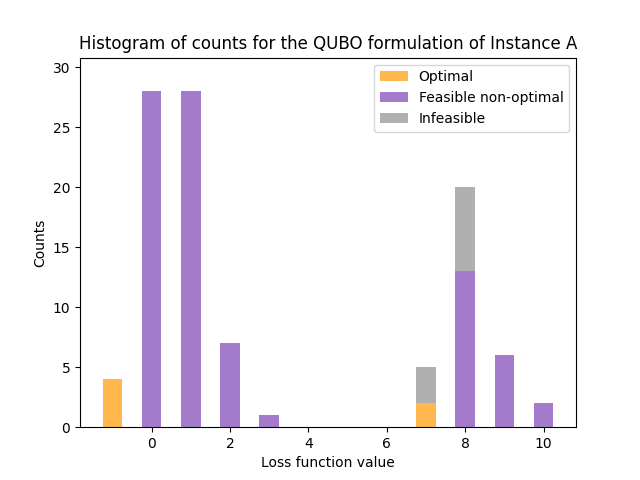}
    }\hspace{0.01\textwidth} 
    \subfigure{
        \includegraphics[width=0.46\textwidth]{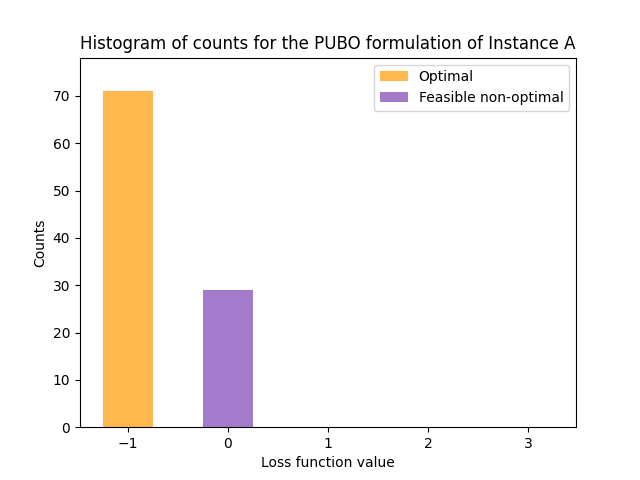}
    }
    \subfigure{
        \includegraphics[width=0.46\textwidth]{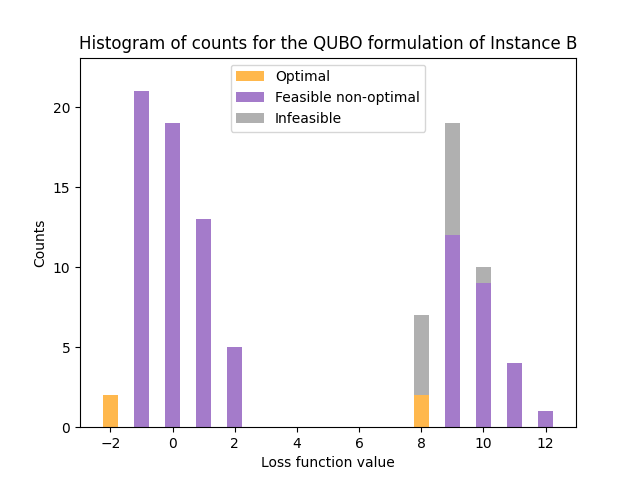}
    }\hspace{0.01\textwidth} 
    \subfigure{
        \includegraphics[width=0.46\textwidth]{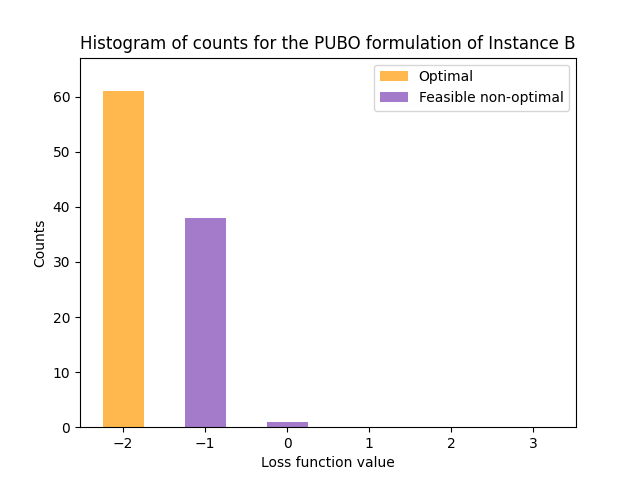}
    }
    \subfigure{
        \includegraphics[width=0.46\textwidth]{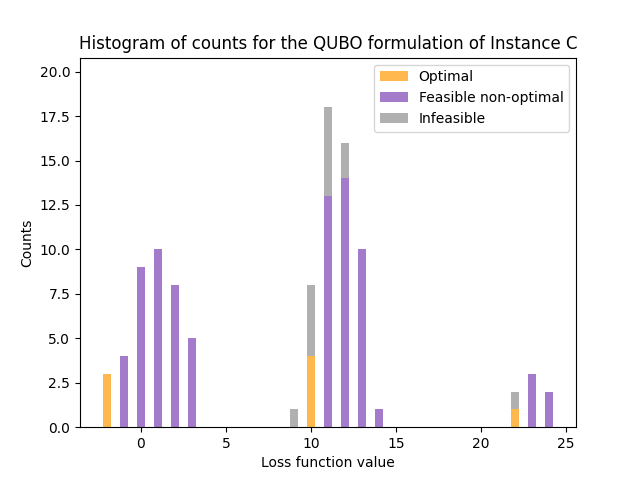}
    }\hspace{0.01\textwidth} 
    \subfigure{
        \includegraphics[width=0.46\textwidth]{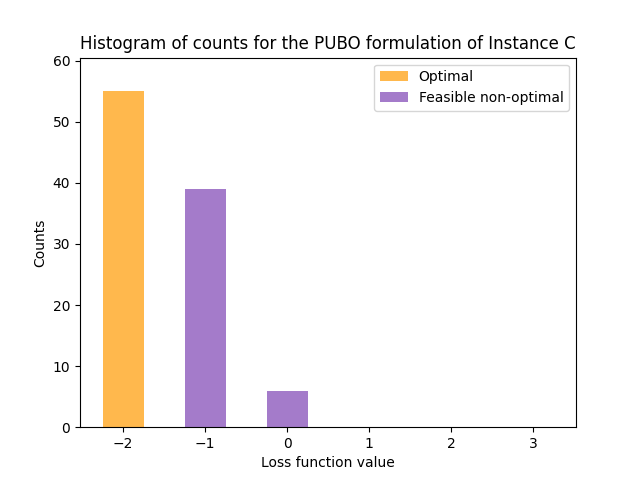}
    }
    \caption{Details of the states returned by running QAOA 100 times for each of the 3 instances (Instance A: row 1, Instance B: row 2, Instance C: row 3) and for each formulation (QUBO: left column, PUBO: right column). The loss function values relate to the unconstrained problem whereas the optimal/feasible non-optimal/infeasible labels' solutions relate to the Extend Bin Packing problem (not yet reformulated).}
    \label{fig:results_OptimBilan_details}
\end{figure}

In Table~\ref{tab:results}, we provide the statistics of the two reformulations (QUBO and PUBO) when solving 100 times each of the three instances with QAOA. The details of the results for the 100 runs are displayed in Figure~\ref{fig:results_OptimBilan_details}. For a run, QAOA returns the best solution encountered during the hybrid optimization, that is the basis state with the smallest loss function value among all the $N_\text{shots}\times N_\text{optim}$ basis states measured in total, where $N_\text{optim}$ is the number of iterations of the classical optimizer, and $N_\text{shots}$ the number of shots to sample the quantum state at each iteration. Here, $N_\text{shots}$ is set to 10 to maintain a reasonable ratio between the sample size and the search space size, thus preserving the relevance of our experiments' scalability. The classical optimizer is chosen to be COBYLA, and no restriction is put on $N_\text{optim}$. Finally, the depth of QAOA is set to 1 to ensure that the circuit depth remains manageable relative to the instance size. For each run, we return the loss function value of the unconstrained problem for the solution given by QAOA, and apply the solution to the initial Extended Bin Packing problem to return if it is optimal, feasible non-optimal, or infeasible.

\begin{table}[!h]
\centering
\small
\begin{tabular}{|c|c|c|c|c|c|c|}\hline
        \multirow{2}{*}{} & \multicolumn{2}{c|}{\textbf{Instance A}} & \multicolumn{2}{c|}{\textbf{Instance B}} & \multicolumn{2}{c|}{\textbf{Instance C}} \\ 
        & QUBO & PUBO & QUBO & PUBO & QUBO & PUBO \\ \hline
        \makecell{Number of qubits} & 15 & 7  & 17 & 9 & 20 & 11 \\ \hline
        \makecell{Proportion of\\ optimal solutions} & 6\% & 71\%  & 4\% & 61\% & 8\% & 55\% \\ \hline
        \makecell{Proportion of feasible\\ non-optimal solutions} & 84\% & 29\%  & 84\% & 39\% & 79\% & 45\%\\ \hline
        \makecell{Proportion of\\ infeasible solutions} & 10\% & 0\%  & 13\% & 0\% & 13\% & 0\% \\ \hline
\end{tabular}
\caption{QAOA solution returned for QUBO and PUBO reformulations. The optimality, feasibility non-optimality and infeasibility labels relate to the Extended Bin Packing problem (not yet reformulated).}
\label{tab:results}
\end{table}

Notice that we implement QAOA solely using the pre-coded COBYLA optimizer and building the quantum circuits from Qiskit library to be executed on the classical 32-qubit simulator of IBM (backend named \emph{qasm\_simulator}). We do not use the pre-coded QAOA function provided by Qiskit and implement, for the PUBO formulation, the quantum circuit with the decomposition proposed by~\cite{Grange_VQA} (Proposition 37). The choice of executing the quantum circuits on a simulator instead on a real quantum hardware is motivated by the fact that today, available quantum backends are still noisy. Thus, to exhibit the impact of the formulation of the input problem on QAOA performances, we decide to use a simulator to \emph{isolate} this parameter from noise and evaluate its real impact.

The PUBO formulation gives better results than the QUBO one on the three instances considered here. Table~\ref{tab:results} shows that over 100 runs over each of the 3 instances, QAOA has never returned an unfeasible solution with the PUBO formulation, whereas between 10\% and 13\% of the solutions returned with the QUBO formulation are infeasible depending on the instances. This means that for some runs of the QUBO formulation, none of the feasible solutions were measured during the whole simulation. Moreover, the PUBO formulation gives far more optimal results than the QUBO formulation: depending on the instances, between 55\% and 71\% of the runs with the PUBO formulation gave an optimal solution, while is was only 4\% to 8\% for the QUBO formulation. Figure~\ref{fig:results_OptimBilan_details} provides a first insight into the reasons of those differences. First, notice that the \textit{packets} that can be observed for the QUBO formulation correspond to different values of the penalty functions: each time a penalty function is increased by 1, it adds $\lambda_c$ or $\lambda_u$ (depending on the constraint) to the loss function of the unconstrained problem. The PUBO reformulation should present packets as well, but they are not visible in our tests because the PUBO problem is well optimized, and only solutions belonging to the first packet are returned. Those diagrams show that for the PUBO formulation, the optimization always encounters one of the first best solutions in the feasible space, while for the QUBO formulation, the optimizer seems to get lost in the multiple packets of the constraints, and does not always reach the feasible space. Notice that for the QUBO reformulation, optimal solutions from the point of view of the nominal problem do not necessarily minimize the loss function of the unconstrained problem (as mentioned in Subsection~\ref{subsec:generic_QUBO_transformation} as \emph{gathering} two minimization problems) because additional variables are not always well optimized. This is why feasible solutions not only appear in the first packet, meaning the one of lower cost. For the PUBO formulation, the packets of higher cost contain only infeasible solutions.

\subsection{Discussion} 

The results of the previous subsection show that the PUBO formulation outperforms the QUBO one for the problem at hand. The possible reasons for this behavior are analyzed below.

\paragraph{Impact of the number of qubits.} The first advantage of PUBO reformulations is that they necessitate less qubits than QUBO reformulations. This is a general difference between QUBO and PUBO reformulations, and is not specific to the two particular reformulations which are considered here. Indeed, in the general case, QUBO reformulations use additional variables, thus additional qubits, to integrate the constraints (see Section~\ref{sec:generic_reformulation}). This difference has two impacts on the quantum part of the hybrid optimization of QAOA. First, PUBO reformulations would be possible on smaller-scale quantum devices. Second, the quantum space to sample is much smaller for PUBO than QUBO as its size evolves exponentially with the number of qubits. Thus, with the same sample size (10 shots here), the probability to measure a good or even optimal solution during the optimization process is enhanced for PUBO reformulations. In our experiment, it is undeniable that the number of qubits plays a role in the advantage of PUBO over QUBO, facilitating the measurement of good solutions at each iteration of the classical optimizer. Our QUBO reformulation uses roughly $(m + n\cdot \log_2(\CMax))$ more qubits compared to the PUBO one. For the three instances of our experiment, it correspond approximately to doubling the number of qubits between the PUBO and the QUBO formulations.

\paragraph{Impact of additional variables.} As mentioned above, using additional variables for QUBO reformulations has a direct impact on the number of qubits and thus on the size of the quantum circuit. However, it also impacts the classical optimization process. In the nominal constrained optimization problem, the set of constraints defines the feasible space by drawing boundaries in the optimization space. In the PUBO reformulation, the same optimization space is kept, but an important cost is assigned to the loss function when crossing the constraints boundaries. On the other side, the QUBO reformulation increases the dimension of the optimization space with additional variables to assign this cost to the loss function. Consequently, the problem has to be optimized over those new variables as well. When those variables are not correctly optimized, they take the advantage over the cost of the nominal problem in the new loss function because they result in huge penalty terms. Consequently, the relevant information, meaning the value of the nominal loss function, gets \textit{blured}. This phenomenon can be observed in our experiments: in Figure~\ref{fig:results_OptimBilan_details} the optimal solution of the initial problem (in orange) can have a loss function value different from the smallest one possible and even larger than other non-optimal solutions for the QUBO reformulation. This is due to the introduction of the additional variables $(s,r^\text{bin})$ in~\eqref{problem:SCP_extended_QUBO} that can take values such that the penalty terms are non-zero, even if the decision variables $(x,y)$ representing a solution of the initial problem are optimal. Thus, finding an optimal solution for the initial problem embedded in a non-optimal solution of the unconstrained problem is a \emph{coincidence} and has nothing to do with optimization. It points out a limit of the QUBO formulation with additional variables that \emph{blur} the loss function. Notice that this phenomenon does not appear for the PUBO formulation~\eqref{problem:SCP_extended_PUBO} (on the right column) because there are no additional variables.

\paragraph{Impact of the optimization landscape.} As explained before, when a constrained problem is reformulated into an unconstrained one, the uncrossable boundaries of the constraints are replaced by a supplementary cost in the loss function. With the QUBO reformulation, this cost can take a large range of values, and the loss function has to be optimized over it even if it is not meaningful for the problem at stake. As long as the feasible space is not reached, the cost that is optimized has no meaning relatively to the nominal problem. The PUBO problem proposed here is meant to mimic the initial constrained problem, by adding an almost constant step to the optimization landscape when crossing the feasible space boundaries. To achieve this property, the penalty functions are built to take only two values: 0 if the associated constraint is satisfied, 1 otherwise. Depending on the classical optimizer, this simplified landscape could facilitate the optimization. In our experiments, it could be partly responsible for the better results of the PUBO formulation, because COBYLA does not get lost in optimizing over the infeasible space as for the QUBO formulation.
\newline

To summarize, several factors explain the advantage of our PUBO formulation with binary constraints over the generic QUBO formulation. First, the QUBO formulation requires more qubits than the PUBO formulation for the same instance of the nominal problem, leading to a lower probability to measure a good basis state during the optimization. Second, the additional variables necessary for the QUBO reformulation have to be optimized as well. They \textit{blur} the relevant loss function, and add dimensions to the space to be optimized classically. Eventually, these new variables also give rise to a large range of values of the penalty terms, leading to an optimization landscape of the QUBO reformulation which does not reflect the one of the unconstrained initial problem. Those two last points result in an optimization landscape that the classical optimizer does not seem to handle very well according to the three instances solved above. Thus, solving the Extended Bin Packing problem with QAOA achieves better performances with the unconstrained formulation~\eqref{problem:SCP_extended_PUBO} than with the~\eqref{problem:SCP_extended_QUBO} formulation on these instances.

\section{Conclusion}
\label{sec:conclusion}

In this paper, we present generic reformulations of unconstrained combinatorial problems to be solved by Variational Quantum Algorithms, and among them, QAOA. Specifically, we introduce a Polynomial Unconstrained Binary Optimization (PUBO) formulation and two generic Quadratic Unconstrained Binary Optimization (QUBO) formulations, valid for every constrained polynomial optimization problem. In addition, we take advantage of the large possibilities offered by PUBO formulations to propose a formulation of optimization problems with linear constraints into a PUBO problem with binary penalty values. We apply this last formulation and the most common QUBO one on a simplified railway timetabling problem stemming from an industrial problem, and solve each of them with QAOA on several small instances.
Our numerical experiments show that the PUBO performs better, drawing a trend for larger instances which haven't been solved yet by lack of high-quality quantum resources. The results speak in favor of the PUBO formulation because it requires less qubits than the QUBO formulation, do not introduce additional variables and seems to present an optimization landscape nicer for the classical optimizer. As soon as large and high-quality quantum hardwares are available, we will be able to deal with the nominal (and more complex) railway timetabling problem of SNCF, comparing PUBO and QUBO formulations on it. This would be the first step of considering quantum algorithms to solve industrial problems. 

However, the two reformulations we compare in this paper involve many simultaneous effects, and further studies should be necessary to clearly distinguish their respective impacts. On the one hand, we should study the impact of additional variables by comparing a generic PUBO formulation (Subsubsection~\ref{subsubsec:generic_PUBO_transformation_general}) and a generic QUBO formulation (Subsection~\ref{subsec:generic_QUBO_transformation}). In order to distinguish the impact of additional variables on the classical optimization from the impact of the number of qubits, QAOA could be run with an exact classical computation of the expectation value of the loss function at each iteration of the optimizer. On the other hand, we should analyze the effect of the shape of the penalty functions and the range of values they can take (binary or not) by comparing the two PUBO formulations presented in Subsection~\ref{subsec:generic_PUBO_transformation}.
Future work should also be dedicated to assess the impact of the noise by implementing QAOA on current quantum computers. Working on noisy quantum computers might be advantageous for the QUBO formulation because of a simpler quantum circuit with a smaller depth.


\paragraph{Acknowledgments.} This work has been partially financed by the ANRT through the PhD number 2021/0281 with CIFRE funds.

\bibliographystyle{apalike}
\bibliography{ref}

\appendix

\section{Omitted proofs}

\subsection{Proof of Property~\ref{prop:generic_cst_equality}}
\label{appendix:proof_cst_equality}
\GenericCstEquality*

\proof
Let $x\in\{0,1\}^n$. 
Let us begin with the case $c=0$. 
\begin{itemize}
    \item If $x$ satisfies~\eqref{eq:constraint_eq}, then $\sum_{i=1}^n x_i = 0$ by definition, and it directly results that $\piEq_0(x) = 0$. 
    \item If $x$ violates~\eqref{eq:constraint_eq}, then we note $\alpha := \sum_{i=1}^n x_i$. By definition of the violation, $\alpha\in\llbracket 1,n\rrbracket$. Thus, because for $k>\alpha$, any product of $k$ variables is equal to 0,
\begin{align*}
    \piEq_0(x) &= \sum_{k=1}^\alpha (-1)^{k+1}\sum_{i \in \EnsBinomial{\alpha}{k}} x_{i_1}\ldots x_{i_k}\\
    &= \sum_{k=1}^\alpha (-1)^{k+1}\binom{\alpha}{k}\\
    &= -\left( \sum_{k=0}^\alpha (-1)^{k}\binom{\alpha}{k} -1 \right)\\
    &= 1 - \sum_{k=0}^\alpha (-1)^{k} 1^{\alpha - k}\binom{\alpha}{k} = 1\,.~~\text{ (Newton binomial formula)}
\end{align*}
\end{itemize}

Let us next consider the case $c\in \N^*$. 
\begin{itemize}
    \item If $x$ satisfies~\eqref{eq:constraint_eq}, then it exists $c$ variables equal to 1. Let us refer to them as $\hat{x}_1,\ldots,\hat{x}_c$. Thus, $\sum_{i\in\EnsBinomial{n}{c}} x_{i_1}\ldots x_{i_c} = \hat{x}_{1}\ldots \hat{x}_{c} = 1$, and $\sum_{i\in\EnsBinomial{n}{k}} x_{i_1}\ldots x_{i_k} = 0$ for $k>c$. It results that 
    \begin{equation*}
        \piEq_c(x) = 1 + (-1)^{c-c+1}\binom{c}{c}\cdot 1 = 1-1 = 0\,.
    \end{equation*}
    \item If $x$ violates~\eqref{eq:constraint_eq}, then by definition $\sum_{i=1}^n = \alpha \neq c$:
    \begin{itemize}
        \item If $\alpha < c$, thus $\sum_{i\in\EnsBinomial{n}{k}} x_{i_1}\ldots x_{i_k} = 0$ for any $k\geq c$, leading to $\piEq_c(x) = 1 - 0 = 1$. 
        \item If $\alpha > c$, thus $\sum_{i\in\EnsBinomial{n}{k}} x_{i_1}\ldots x_{i_k} = \binom{\alpha}{k}$ for any $c\leq k\leq \alpha$, and $\sum_{i\in\EnsBinomial{n}{k}} x_{i_1}\ldots x_{i_k} = 0$ for any $k>\alpha$. It results that 
        \begin{equation*}
            \piEq_c(x) = 1 - \sum_{k = c}^\alpha (-1)^{k-c+1}\binom{k}{c}\binom{\alpha}{k}\,,
        \end{equation*}
        where we can show by manipulating factorials that \begin{equation*}
            \binom{k}{c}\binom{\alpha}{k} = \binom{\alpha}{c}\binom{\alpha - c}{k - c}\,.
        \end{equation*}
        Thus, \begin{align*}
            \piEq_c(x) &= 1 - \sum_{k = c}^\alpha (-1)^{k-c+1}\binom{\alpha}{c}\binom{\alpha - c}{k - c}\\
            &= 1 + \binom{\alpha}{c}\sum_{k = 0}^{\alpha-c} (-1)^{k}\binom{\alpha - c}{k} \\
            &= 1 + \binom{\alpha}{c} (1 - 1)^{\alpha - c} = 1\,.
        \end{align*}
    \end{itemize}
\end{itemize}
\endproof

\subsection{Proof of Property~\ref{prop:generic_cst_inequality_inf}}
\label{appendix:proof_cst_inequality_inf}
\GenericCstInequalityInf*

\proof 
Let $x\in\{0,1\}^n$.
\begin{itemize}
    \item If $x$ satisfies~\eqref{eq:constraint_inf}, then $\sum_{i\in\EnsBinomial{n}{k}} x_{i_1}\ldots x_{i_1} = 0$ for any $k>c$. Thus, $\piInf_c(x) = 0$.
    \item If $x$ violates~\eqref{eq:constraint_inf}, let us note $\sum_{i=1}^n x_i = \alpha > c$. Thus, 
    \begin{equation*}
        \sum_{i\in\EnsBinomial{n}{k}} x_{i_1}\ldots x_{i_k} = \begin{cases}
            \binom{\alpha}{k}~~\text{ for } k\in\llbracket c+1, \alpha\rrbracket\\
            0~~\text{ for } k>\alpha
        \end{cases}
    \end{equation*}
    It results that $$\piInf_c(x) = \sum_{c+1}^\alpha (-1)^{k-c+1}\binom{k-1}{c}\binom{\alpha}{k}\,.$$ Next, we prove by recurrence over $c\in\llbracket 0, \alpha-1\rrbracket$ the proposition $$\mathcal{P}(c) : \text{``}\piInf_c(x) = 1, \text{ for all } x\in\{0,1\}^n \text{ violating}~\eqref{eq:constraint_inf}\text{''}\,.$$

    \textit{Initialization:} $\mathcal{P}(\alpha - 1)$ is True. Indeed, for $x \text{ violating}~\eqref{eq:constraint_inf}$, 
    
    $$\piInf_{\alpha - 1}(x) = (-1)^{\alpha-(\alpha - 1) +1} \binom{\alpha-1}{\alpha-1}\binom{\alpha}{\alpha} = 1\,.$$ 

    \textit{Recurrence:} Let $c\in\llbracket1,\alpha-1\rrbracket$, and let us assume that $\mathcal{P}(c)$ is True. Let us show that $\mathcal{P}(c-1)$ is also True. For $x \text{ violating}~\eqref{eq:constraint_inf}$,
    \begin{align*}
        \piInf_{c- 1}(x) &= \sum_{k=c}^\alpha (-1)^{k-c} \binom{k-1}{c-1}\binom{\alpha}{k}\\
        &= (-1)^{c-c}\binom{c-1}{c-1}\binom{\alpha}{c} + \sum_{k=c+1}^\alpha (-1)^{k-c} \binom{k-1}{c-1}\binom{\alpha}{k}\\
        &= \binom{\alpha}{c} - \sum_{k=c+1}^\alpha (-1)^{k-c +1} \left(\binom{k}{c} - \binom{k-1}{c}\right)\binom{\alpha}{k}~~\text{(Pascal's triangle)}\\
        &= \binom{\alpha}{c} - \sum_{k=c+1}^\alpha (-1)^{k-c +1}\binom{k}{c}\binom{\alpha}{k} + \underbrace{\sum_{k=c+1}^\alpha (-1)^{k-c +1}\binom{k-1}{c}\binom{\alpha}{k}}_{\piInf_c(x)}\\
        &= \binom{\alpha}{c} - \sum_{k=c+1}^\alpha (-1)^{k-c +1}\binom{k}{c}\binom{\alpha}{k} + 1\,.~~\text{(Recurrence hypothesis)}\\
    \end{align*}

    Identically to the proof of the penalty term of constraint~\eqref{eq:constraint_eq} for the case $c\in\N^*$, we can show that $$\sum_{k=c+1}^\alpha (-1)^{k-c +1}\binom{k}{c}\binom{\alpha}{k} = \binom{\alpha}{c}\,.$$ Thus, it results that $\piInf_{c- 1}(x) = \binom{\alpha}{c} - \binom{\alpha}{c} + 1 = 1$.
    
\end{itemize}
\endproof

\subsection{Proof of Property~\ref{prop:generic_cst_inequality_sup}}
\label{appendix:proof_cst_inequality_sup}
\GenericCstInequalitySup*

\proof
Let $x\in\{0,1\}^n$. It is sufficient to show that $\piSup_c(x) = \piEq_c(x) - \piInf_c(x)$, for $c\in\N^*$. Indeed, if we note $\alpha := \sum_{i=1}^n x_i$, we have the following results. 
\begin{itemize}
    \item If $x \text{ satisfies}~\eqref{eq:constraint_sup}$ (meaning that $\alpha \geq c)$:
    \begin{itemize}
        \item If $\alpha = c$\,: $\piEq_c(x) = \piInf_c(x) = 0$, thus $\piSup_c(x) = 0$.
        \item If $\alpha >c\,$: $\piEq_c(x) = \piInf_c(x) = 1$, thus $\piSup_c(x) = 0$.
    \end{itemize}
    \item If $x \text{ violates}~\eqref{eq:constraint_sup}$ (meaning that $\alpha<c)$, we have $\piEq_c(x) = 1$ and $\piInf_c(x) = 0$, leading to $\piSup_c(x) = 1$.
\end{itemize}
Thus, it remains to prove that $\piSup_c(x) = \piEq_c(x) - \piInf_c(x)$.
\begin{align*}
    \piEq_c(x) - \piInf_c(x) &= 1 + \sum_{k=c}^n (-1)^{k-c+1} \binom{k}{c}\sum_{i \in \EnsBinomial{n}{k}} x_{i_1}\ldots x_{i_k} - \sum_{k=c+1}^n (-1)^{k-c+1} \binom{k-1}{c}\sum_{i\in \EnsBinomial{n}{k}} x_{i_1}\ldots x_{i_k}\\
    &= 1 + (-1)^{c-c+1}\binom{c}{c}\sum_{i \in \EnsBinomial{n}{c}} x_{i_1}\ldots x_{i_k} + \sum_{k=c+1}^n (-1)^{k-c+1} \left(\binom{k}{c} - \binom{k-1}{c}\right)\sum_{i\in \EnsBinomial{n}{k}} x_{i_1}\ldots x_{i_k}\\
    &= 1 - (-1)^{c-c+1}\binom{c-1}{c-1}\sum_{i \in \EnsBinomial{n}{c}} x_{i_1}\ldots x_{i_k} + \sum_{k=c+1}^n (-1)^{k-c+1} \binom{k-1}{c-1}\sum_{i\in \EnsBinomial{n}{k}} x_{i_1}\ldots x_{i_k}\\
    &= 1 + \sum_{k=c}^n (-1)^{k-c+1} \binom{k-1}{c-1}\sum_{i\in \EnsBinomial{n}{k}} x_{i_1}\ldots x_{i_k}\\
    &= \piSup_c(x)\,.
\end{align*}
\endproof

\end{document}